\documentclass[12pt]{article}
\pdfoutput=1
\usepackage{fullpage}
\usepackage{amsmath,amssymb,amsthm}
\usepackage{xcolor}
\usepackage{graphics,graphicx}
\usepackage{algorithmic}
\usepackage{algorithm}

\usepackage{hyperref}
\usepackage{setspace}

\begin{document}

\def\E{{\mathbb E}}
\def\P{{\mathbb P}}
\def\R{{\mathbb R}}
\def\Z{{\mathbb Z}}
\def\V{{\mathbb V}}
\def\N{{\mathbb N}}
\def\etatnx{{\frac{\eta_t^N(x)}{N}}}
\def\etanx{{\frac{\eta^N(x)}{N}}}
\def\etany{{\frac{\eta^N(y)}{N}}}
\def\etatny{{\frac{\eta_t^N(y)}{N}}}
\def\ngoes{{N \to \infty}}
\newcommand\cW{\mathcal{W}}
\newcommand{\cS}{\mathcal{S}}
\newcommand{\MP}[1]{\marginpar{\tiny #1}}
\newcommand{\vep}{\varepsilon}
\newcommand{\cA}{\mathcal{A}}
\newcommand{\cY}{\mathcal{Y}}
\newcommand{\cYjk}{\mathcal{Y}_{jk}}
\newcommand{\cOjk}{\mathcal{O}_{jk}}
\newcommand{\cOo}{\mathcal{O}_{00}}
\newcommand{\cYo}{\mathcal{Y}_{00}}
\newcommand{\cP}{\mathcal{P}}
\newcommand{\Lso}{{L}^{small}_{00}}
\newcommand{\Rso}{{R}^{small}_{00}}
\newcommand{\cL}{\mathcal{L}}
\newcommand{\cE}{\mathcal{E}}

\def\ind{\ensuremath{{{1}}\kern-0.29em{{\hbox{\rm I}}}}}
\def\Ljk{{L_{jk}}}
\def\Rjk{R_{jk}}
\def\Ljks{{L^{small}_{jk}}}
\def\Rjks{R^{small}_{jk}}
\def\Ojk{\mathcal{O}_{jk}}
\def\O00{\mathcal{O}_{00}}

\def\Ejk{{\mathcal{E}_{jk}}}
\def\a{{\alpha}}
\def\b{{\beta}}

\def\l{{\lambda}}
\def\qed{\hfill $\blacksquare$}

\def\({{\Bigl(}}
\def\){{\Bigr)}}

\newtheorem{teo}{Theorem}
\newtheorem{lemma}[teo]{Lemma}
\newtheorem{propo}[teo]{Proposition}
\newtheorem{conj}{Conjecture}
\newtheorem{alg}{Algorithm}[subsection]
\newtheorem{cor}[teo]{Corollary}
\newtheorem{de}{Definition}
\newtheorem{rem}[teo]{Remark}

\title{Contact process in a wedge}
\author{J. Theodore Cox\thanks{ Supported in part by NSF Grant
  No. 0803517}\\ {\it Syracuse University} 
\and
 Nevena Mari\'{c}\thanks{ Supported in part by NSF Grant
  No. 0803517}\\ {\it University of Missouri - St. Louis}
 \and
  Rinaldo Schinazi\thanks{ Supported in part by NSF Grant
  No. 0701396}\\ {\it University of Colorado - Colorado Springs}}

\maketitle

\noindent
\abstract We prove that the supercritical
one-dimensional contact process survives in certain wedge-like
space-time regions, and that when it
survives it couples with the unrestricted contact process
started from its upper invariant measure. As an application we show
that a type of weak coexistence is possible in the nearest-neighbor
``grass-bushes-trees'' successional model introduced in \cite{swindle}.
\endabstract


\section{Introduction}


The contact process of Harris (introduced in \cite{Ha}) is a well
known model of infection spread by contact.  The
one-dimensional model is a continuous time Markov process
$\xi_t$ on $\{0,1\}^\Z$.  For $x\in\Z$, $\xi_t(x)=1$ means
the individual at site $x$ is infected at time $t$ while
$\xi_t(x)=0$ means the individual is healthy.  Infected
individuals recover from their infection after an
exponential time with mean 1, independently of everything
else. Healthy individuals become infected at a rate
proportional to the number of infected neighbors.
Alternatively, individuals (1's) die at rate one and give
birth onto neighboring empty sites (0's) at rate $\lambda$.
If we let $n_i(x,\xi)=\sum_{y:|y-x|=1} 1\{\xi(y)=i\}$, and
$\l\ge0$ the infection parameter, then the transitions at
$x$ in state $\xi$ are
\begin{equation}\label{stdrates}
1 \to 0 \text{ at rate } 1 \quad\text{ and }\quad
0\to 1  \text{ at rate } \l n_1(x,\xi)\,.
\end{equation}
When convenient we will  identify $\xi\in\{0,1\}^{\Z}$ with $\{x:
\xi(x)=1\}$, and use the notation  $\|\xi\|_i=\sum_x 1\{\xi(x)=i\}$. 

Let $\xi^0_t$ denote the contact process with initial
state $\xi^0_0=\{0\}$. The critical value $\l_c$ is defined by
\begin{equation}
\l_c= \inf\{ \l \geq 0: P(\xi^0_t\ne\emptyset \text{ for all }t\ge
0)>0 \}\,. 
\end{equation} 
It is well known that $0<\l_c<\infty$, and that in the
supercritical case $\l>\l_c$ there is a unique stationary 
distribution $\nu$
for $\xi_t$, called the upper invariant measure, with the
property 
\[
\nu\bigl(\xi:  \|\xi\|_1 = \infty\bigr) = 1 \,.
\] 
There are also well-defined ``edge speeds.'' Let
$\xi_0^-$ ($\xi_0^+$) be the initial
state given by 
$\xi_0^-=\Z^-$ ($\xi_0^+=\Z^+$), and define the edge
processes 
\begin{equation}
r_t = \max\{x: \xi_t^- (x)=1\}\text{ and }
l_t = \min\{x: \xi_t^+ (x)=1\}\,.
\end{equation}
There is a strictly increasing function
$\alpha:(\l_c,\infty)\to(0,\infty)$  such that 
for $\l>\l_c$ 
\begin{equation}\label{edgespeed}
\lim_{t \to \infty} \frac{r_t}{t}= \a(\l) 
\text{ and }
\lim_{t \to \infty} \frac{l_t}{t}= -\a(\l) ~ ~ a.s. 
\end{equation}
All of the above facts are contained in 
Chapter VI of \cite{bible} and Part I of \cite{bible2}.

We are interested in contact processes for which the infection
is restricted
to certain space-time regions.  For $\cW\subset
\Z\times[0,\infty)$ define the $\cW$-restricted contact
process $\xi^\cW_t$ as follows. First, set
$\xi^\cW_t(x) =0 $ for all $(x,t)\notin \cW$. Second, for
$(x,t)\in \cW$, replace \eqref{stdrates} with 
\begin{equation}\label{wedgerates}
1 \to 0 \text{ at rate } 1 \quad\text{ and }\quad
0\to 1  \text{ at rate } \l\sum_{y:|y-x|=1}
\xi(y)1_\cW(y,t) \,, 
\end{equation}
so that infection spreads only between sites in the wedge.  We will
give an explicit \emph{graphical construction} of $\xi^\cW_t$ in 
Section~2.

For $0<\alpha_l<\alpha_r<\infty$ and $M\ge0$ define the ``wedges'' 
$\cW=\cW(\alpha_l,\alpha_r,M)$ by
\begin{equation}\label{wedge}
\cW = \{ (x,t) \in \Z\times[0,\infty): \alpha_lt\le x\le M+\alpha_rt 
\} \,.
\end{equation}
In view of \eqref{edgespeed}, we will impose the conditions
\begin{equation}\label{wedgespeeds}
\l > \l_c \text{ and }0< \alpha_l<\alpha_r< \alpha(\l)\,.
\end{equation}
Our first result is that survival in wedges is possible.
\begin{teo}\label{survival} Assume \eqref{wedgespeeds} holds,
$\cW = \cW(\alpha_l,\alpha_r,M)$, and $\xi^\cW_0=[0,M]\cap\Z$. Then
\begin{equation}\label{eqn:survival}
\lim_{M\to\infty}P(\xi^\cW_t\ne\emptyset \text{ for all
}t\ge 0 ) = 1 \,.
\end{equation}
\end{teo}

    When $\xi^\cW_t$ survives
it looks like the unrestricted contact process in equilibrium. 
To state this more precisely, let 
\begin{equation}
r^\cW_t = \max\{x: \xi^\cW_t (x)=1\} \text{ and }
l^\cW_t = \min\{x: \xi_t^\cW (x)=1\}\,,
\end{equation}
and let $\xi^\nu_t$ denote the contact process
started in its upper invariant measure $\nu$. 

\begin{teo}\label{equilibrium} Assume \eqref{wedgespeeds},
$\cW =
\cW(\alpha_l,\alpha_r,M)$, and $\xi^\cW_0=
[0,M]\cap \Z$. On the event $\{\xi_t^\cW\ne\emptyset
\text{ for all }t\ge 0\}$, 
\begin{equation}\label{wedgegrowth}
\lim_{t \to \infty} \frac{r^\cW_t}{t}= \alpha_r
\text{ and }
\lim_{t \to \infty} \frac{l^\cW_t}{t}= \alpha_l ~ a.s.
\end{equation}
Furthermore, $\xi^\cW_t$ and $\xi^\nu_t$ can be coupled so
that on the event $\{\xi^{\cW}_t\ne \emptyset \text{ for all
}t\ge 0\}$, 
\begin{equation}\label{nu-couple}
\xi^\cW_t(x) = \xi^\nu_t(x) \text{ for all }
x\in[l^\cW_t,r^\cW_t] \text{ for all large   }t ~ a.s.
\end{equation}
\end{teo}

\begin{rem}\label{lingrowth} By standard arguments using exponential estimates,
$|\xi^\nu_t\cap [at,bt]|\to\infty$
  as $t\to\infty$ with probability one for any $a<b$
(see Theorem~VI.3.33 in \cite{bible}). 
Therefore Theorem~\ref{equilibrium} implies that when 
$\xi^{\cW}_t$ survives, $|\xi^{\cW}_t|\to\infty$ a.s.
\end{rem}

Theorem~\ref{survival} can be used to obtain  information about the
``grass-bushes-trees'' model (GBT) of \cite{swindle}.  In this
model sites are either empty (0), occupied by a bush
(1) or occupied by a tree (2). Both 1's and 2's
turn to 0's at rate one. The 2's give birth at rate $\lambda_2$
on top of 1's and 0's.  The 1's give birth at rate $\lambda_1$
on top of 0's only, and hence are at a disadvantage compared
to 2's. The state space for the process is $\{0,1,2\}^\Z$, and
the nearest-neighbor version of the model makes transitions at
$x$ in state $\zeta$
\begin{equation}\label{gbt-rates}
0\to \begin{cases}
1 & \text{ at rate } \l_1 n_1(x,\zeta)\\
2 & \text{ at rate } \l_2n_2(x,\zeta)
\end{cases}\, \quad
1\to \begin{cases}
0 &\text{ at rate } 1 \\
2 &\text{ at rate } \l_2n_2(x,\zeta)
\end{cases}\,\quad
2 \to 0  \text{ at rate } 1\,.
\end{equation}

A natural question to ask is whether or not coexistence of 1's
and 2's is possible.
It was shown in \cite{swindle} that coexistence is possible
for a non-nearest neighbor version of the model and
appropriate $\l_i$, where coexistence meant that $\zeta_t$ had 
a stationary distribution $\mu$ such that
\begin{equation}\label{s-coexist}
\mu\Bigl( \zeta: \|\zeta\|_i=\infty \text{ for }i=1,2 \} 
\Bigr)=1 \,.
\end{equation}
It was also shown in \cite{swindle} that there is no
stationary distribution satisfying \eqref{s-coexist}
in the nearest-neighbor case for \emph{any} choice of the
$\l_i$. Moreover, if there are infinitely many 2's initially then
for each site there is a last time at which a 1 can be
present.  Nevertheless, it is a consequence of
Theorem~\ref{survival} and the construction used in its proof
that a form of weak coexistence is possible, even starting
from a single 1 and infinitely many 2's. 

\begin{cor}\label{gbt} Let $\zeta_t$ be the GBT process with initial state 
  $\zeta_0$, where $\zeta_0(x)=2$ for $x<0$, $\zeta_0(0)=1$
  and $\zeta_0(x)=0$ for $x>0$. For all $\l_c<\l_2<\l_1$,
\begin{equation}\label{gbtsurv}
P\Bigl(\lim_{t\to\infty}\|\zeta_t\|_1 =\infty \Bigl) >0 \,.
\end{equation}
\end{cor}
The 2's spread to the right at
rate $\alpha(\lambda_2)$, ignoring the 1's, while the 1's
try to spread to the
right at the faster rate $\alpha(\lambda_1)$ . The 1's will be
killed by 2's invading from the left, but
Theorem~\ref{survival} shows that they can survive with
positive probability by moving off to the right in the
space-time region free of 2's.

\begin{rem} (1) With a little more work one can use Theorem 2 to
  say more about the set of of 1's in $\zeta_t$ since it
  dominates wedge-restricted contact processes with positive
  probability. (2) Non-oriented percolation in various subsets
  of $\Z^d$ has been studied by others (e.g. see \cite{Grim}
  and \cite{CC}), but as
  far as we are aware our results on oriented percolation are new.
\end{rem}

In Section 2 we give the standard graphical construction due
to Harris, then prove Theorem~\ref{survival} in Section 3,
Theorem~\ref{equilibrium} in Section 4, and
Corollary~\ref{gbt} in Section 5.  

\section{The graphical representation}\label{harris}

For $x\in \Z$ let $\{T_n^x:n\geq 1\}$ be the arrival times of
a Poisson process with rate 1, and for all pairs of
nearest-neighbor sites $x,y$ let $\{B_n^{x,y}:n\geq 1\}$ be
the arrival times of a Poisson process with rate $\l$.  The
Poisson processes $T^x,B^{x,y}$, $x,y\in \Z$, are all
independent.  At the times $T_n^x$ we put a $\delta$ at site
$x$ to indicate a death at $x$, 
and at the times $B_n^{x,y}$ we draw an
arrow from $x$ to $y$, indicating that a 1 at $x$
will give birth to a 1 at $y$. For $0\le s<t$ and sites $x,y$ we
say that there is an active path up from $(x,s)$ to $(y,t)$ if
there is a sequence of times $t_0=s\le t_1<t_2<\dots<t_{n}\le
t_{n+1}=t$ and a sequence of sites $x_0=x,x_1,\dots,x_n=y$
such that
\begin{enumerate}
\item for $i=1,2\dots,n$, $|x_i-x_{i-1}|=1$ and there is an arrow from $x_{i-1}$ to
  $x_i$ at time $t_i$
\item for $i=0\dots,n$, the time segments
  $\{x_i\}\times[t_{i},t_{i+1}]$ do not contain any
  $\delta$'s
\end{enumerate}
By default there is always an active path up from $(y,t)$ to
$(y,t)$.  For a space-time region $\cW\subset
\Z\times[0,\infty)$ we define $\xi^{\cW}_t$, the contact
process restricted to $\cW$, as follows. Given an initial
state $\xi_0\subset\{x:(x,0)\subset \cW\}$, set $\xi_t(y) = 0$
for all $(y,t)\notin \cW$. If there is a site $x$ with
$\xi_0(x)=1$ and an active path up from $(x,0)$ to $(y,t)$
lying entirely in $\cW$ set $\xi_t^{\cW}(y)=1$, otherwise
set $\xi^{\cW}_t(y)=0$.  For $\cW=\Z\times[0,\infty)$ we will
write $\xi_t$ and refer to it as the unrestricted
process.

We may also construct the GBT process $\zeta_t$ with the
above Poisson processes and the help of some additional independent
coin flips.  Fix $\lambda_c<\lambda_2<\lambda_1$, and suppose
$\lambda=\lambda_1$ in the construction just given.
Independently of everything else, label the arrows determined
by the $B^{xy}_n$ with a ``1-only'' sign with probability
$(\lambda_1-\lambda_2)/\lambda_1$. Call an active path up from
$(x,s)$ to $(y,t)$ a 2-path if none of its arrows are 1-only
arrows.  Given $\zeta_0$, we may now construct $\zeta_t$ as
follows. First, for all $t>0$ and $x\in\Z$, put $\zeta_t(x)=2$
if for some site $y$ with $\zeta_0(y)=2$ there is an active
2-path up from $(y,0)$ to $(x,t)$.  Next, for all other
$(x,t)$ put $\zeta_t(x)=1$ if for some site $y$ with
$\zeta_0(y)=1$ there is an active path up from $(y,0)$ to
$(x,t)$ with the property that no vertical segments in the
path contain a point $(z,u)$ such that
$\zeta_u(z)=2$. Otherwise set $\zeta_t(x)=0$.  A little
thought shows that $\zeta_t$ is the GBT process with the rates
given in \eqref{gbt-rates}. The process of 2's is a contact
process with infection parameter $\lambda_2$, and in the
absence of 2's, the process of 1's is a contact process with
infection parameter $\lambda_1$.

\section{Proof of Theorem~\ref{survival}}
\paragraph*{The space-time regions $\cYjk$.}
We will modify somewhat the standard approach of constructing
a mapping from appropriate space-time regions of the
construction just given to an oriented-percolation model 
with the property that survival of the percolation process
implies survival of the contact process. 
We will call the regions $\cYjk$, they will be defined
using the parallelograms of Section~VI.3 of \cite{bible}. 

\begin{figure}[htp]
\begin{center}
\includegraphics[totalheight=0.15\textheight]{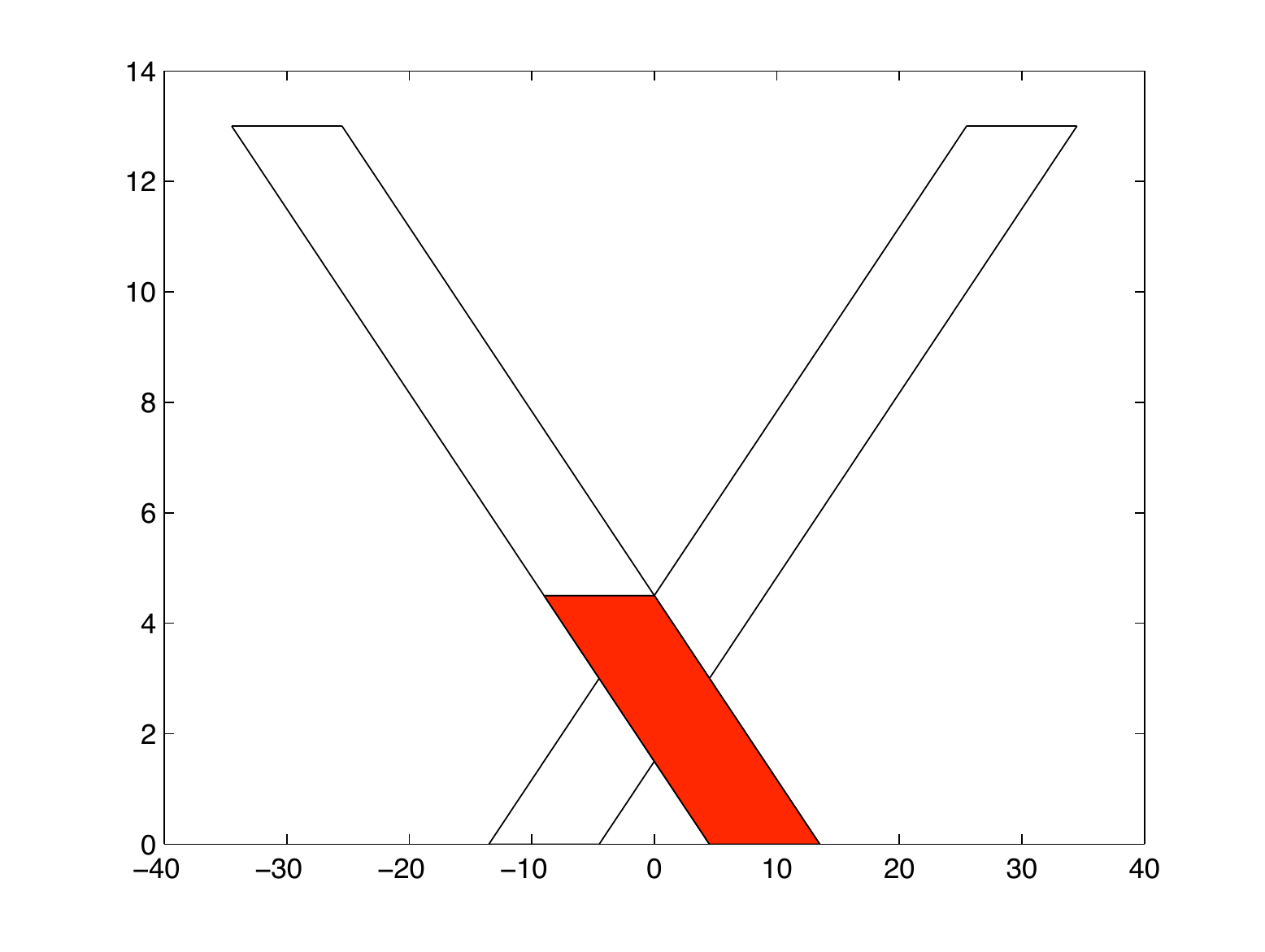}
\caption{Large parallelograms $L_{00} $ and $ R_{00}$. The
  shaded region is $L_{00}^{small}$.} \label{bigsmall} 
\end{center}
\end{figure}

Let $\cL$ be the lattice $\cL=\{ (j,k) \in \Z^2: k \geq 0
~\mbox{and $j+k$ is even}\}$ with norm
$\|(j,k)\|=1/2(|j|+|k|)$.  Fix $0 < \b < \a/3$ and $M>0$ so
that $M \b/2$ and $ M\a$ are integers. Later we will set
$\alpha=\alpha(\lambda)$ and take $\beta$ small. For $(j,k)
\in \cL$, $\Ljk$ and $\Rjk$ are the ``large'' space-time
parallelograms in $\Z \times [0,\infty)$ given by:
\[
\Ljk = M(j(\alpha-\beta),k) + L_{00}, \quad 
\Rjk =  M(j(\alpha-\beta),k) + R_{00}
\]
where 
\begin{align*}
  L_{00} &= \{ (x,t) \in \Z \times [ 0, M(1+ {\b}/{\a})]
 : M\beta/2 \leq x+\a t \leq 3M\b/2 \}\\
R_{00} &= \{ (x,t) \in \Z \times [ 0, M(1+
{\b}/{\a})]: -3M\b/2 \leq x-\a t \leq -M\b/2 \} \;.
\end{align*}
We will also need the ``small'' parallelograms
\[
\Ljks = M(j(\alpha-\beta),k) + L^{small}_{00}, \quad 
\Rjks =  M(j(\alpha-\beta),k) + R^{small}_{00}
\]
where
\begin{align*}
  \Lso &= \{ (x,t) \in \Z \times [ 0, M \frac{3\b}{2\a}]
 : M\beta/2 \leq x+\a t \leq 3M\b/2 \}\\
\Rso &= \{ (x,t) \in \Z \times [ 0, 
M\frac{3\b}{2\a}]: -3M\b/2 \leq x-\a t \leq -M\b/2 \} \;.
\end{align*}
It is important to note that $\Lso\subset L_{00}$,
$\Rso\subset R_{00}$, 
and 
\[
\Rjk \cap \Ljk = \Rjk \cap \Ljks =  \Rjks\cap \Ljk \;,
\]
as shown in Figure \ref{bigsmall}. 
 
\begin{figure}[htp]
\hfill\includegraphics[totalheight=0.15\textheight]{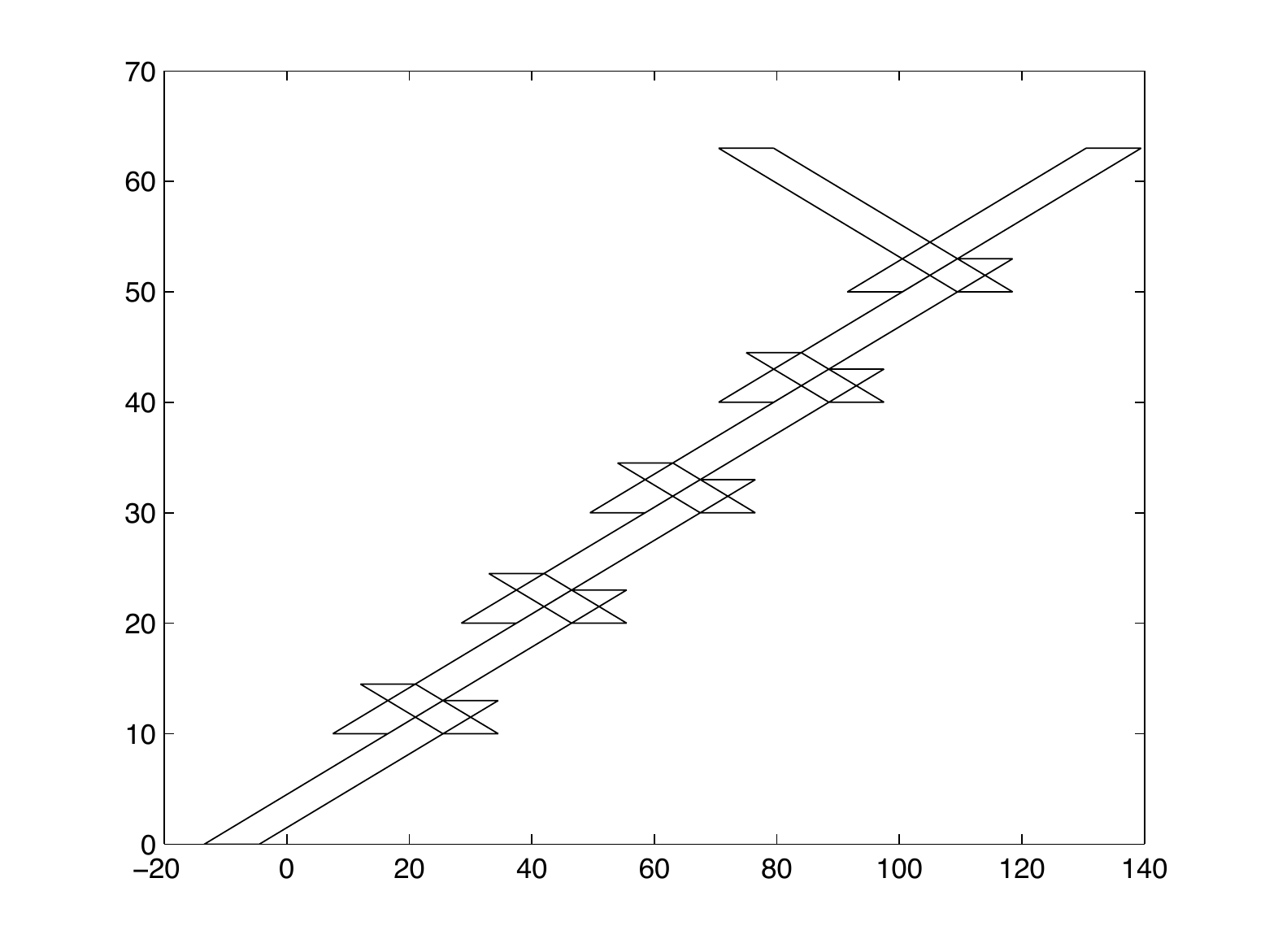}
\hfill
\includegraphics[totalheight=0.15\textheight]{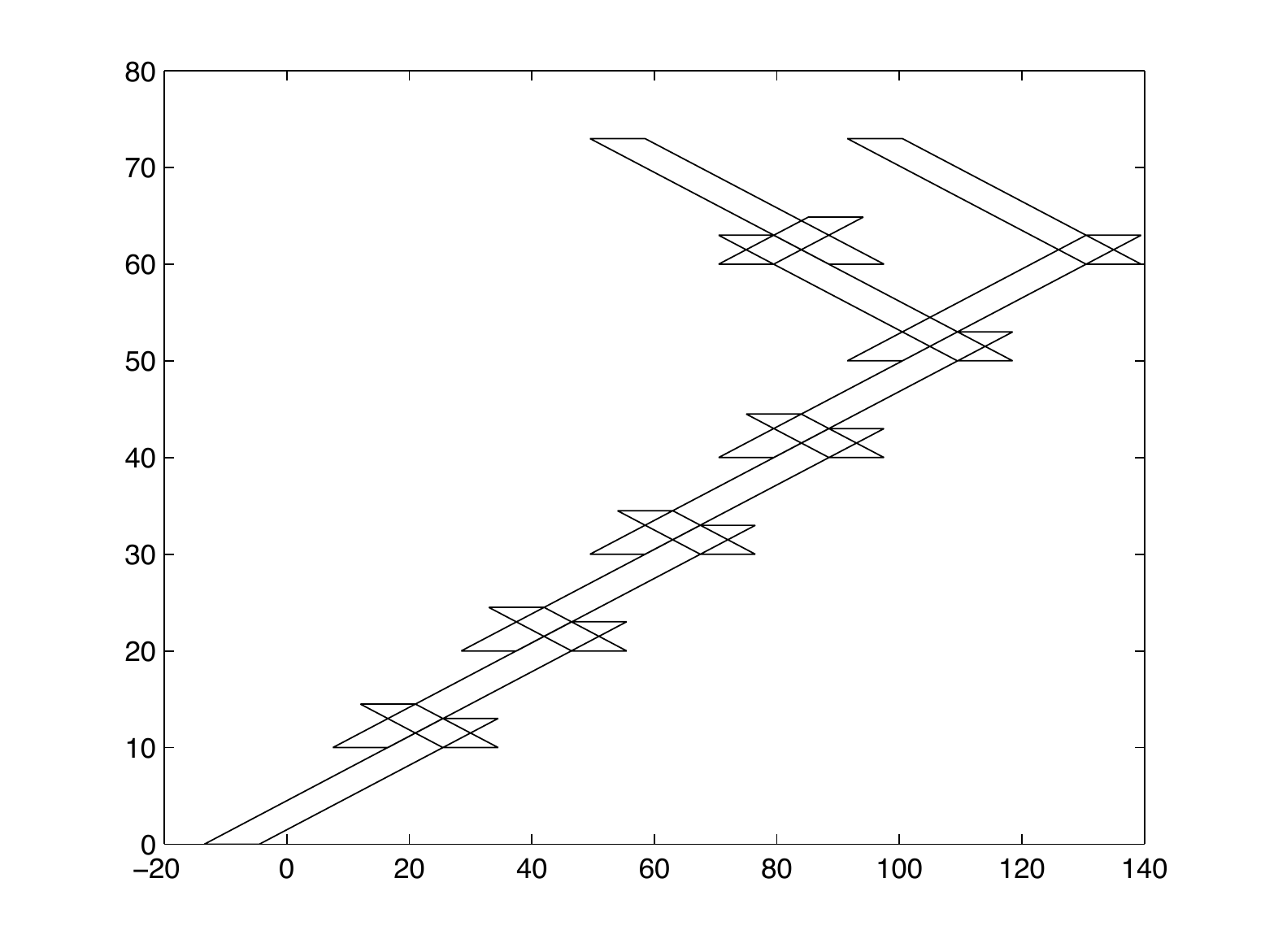}
\hfill
\includegraphics[totalheight=0.15\textheight]{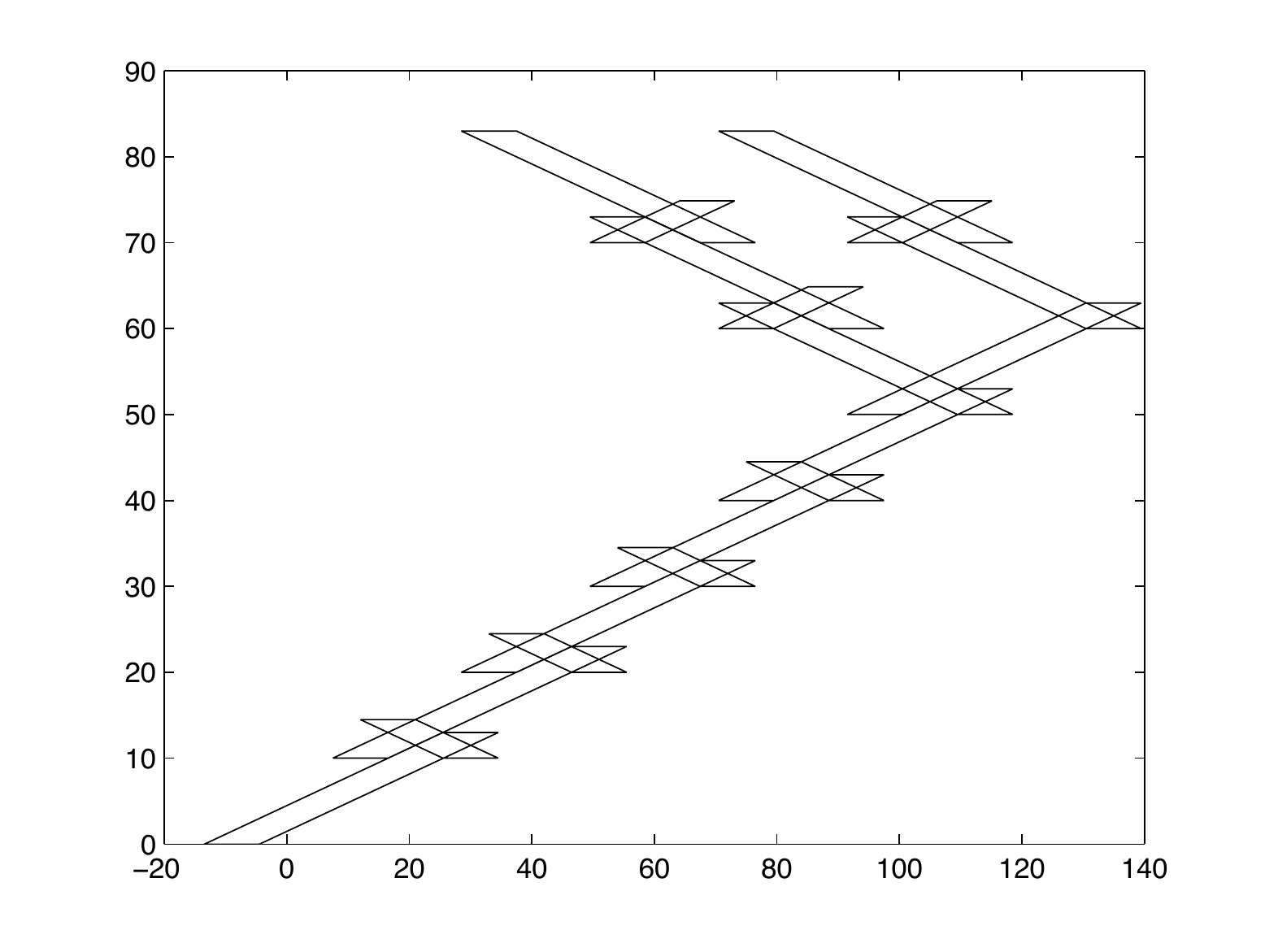}
\hfill
\caption{$\cYo$  with $\ell=5$, $d=0,1,2$.}\label{obj}
\end{figure}

We can now define the new objects $\cYjk$ which will be used
to construct our oriented percolation process. As is the
case with the parallelograms, the $\cYjk$ will be certain translates of
$\cYo$, and depend on two fixed integers $\ell,d$
which satisfy $\ell\ge 2$ and $d\ge 0$ with $\ell>d$.
We will form $\cYo$ by sticking together $\ell$ big right
parallelograms, connected with appropriate small left
parallelograms, and then two branches of $d$ and $d+1$ big
left parallelograms connected by small right
parallelograms. Figure \ref{obj} shows examples of $\cYo$
with parameters $\ell=5$ and $d=0,1,2$.
It seems simplest to define $\cYo$ in
stages, beginning with 
$\cYo^0=R_{00}$.
\begin{enumerate}
\item Attach $\ell$ big right parallelograms with $\ell$
  small parallelograms to connect them: 
\[
\cYo^1 = \cYo^0\cup \Bigl(\bigcup_{i=1}^\ell ( R_{ii}  \cup
L^{small}_{ii})\Bigr) \,.
\]
\item Attach one big left parallelogram:
 $\cYo^2  = \cYo^1 \cup L_{\ell,\ell}$.

\item If $d=0$ set $\cYo=\cYo^2$. If $d\ge 1$, attach another big left
  parallelogram:
\[
\cYo^3 = \cYo^2 \cup  L_{\ell+1,\ell+1}\;.
\]
\item If $d=1$, attach another big left and small right
  parallelogram:
\[
\cYo^4 = \cYo^3 \cup ( L_{\ell-1,\ell+1}\cup
R^{small}_{\ell-1,\ell+1} )
\]
and set $\cYo = \cYo^4$. If {$d \geq 2$}, attach two
branches, to reach ``height'' $\ell+d+1$, of big left
parallelograms with small right parallelograms as
connectors:
\[
\cYo^4 = \cYo^3 \cup \Bigl(
\bigcup_{i=0}^{d-1} 
(L_{\ell-i,\ell+i} \cup R^{small}_{\ell-i,\ell+i} )\cup
(L_{\ell+1-i,\ell+1+i} \cup R^{small}_{\ell+1-i,\ell+1+i})
\Bigr) \;.
\]
\item If $d\ge 2$,  attach a final big left parallelogram
  and small right parallelogram:
\[
\cYo^5= \cYo^4 \cup L_{\ell-d,\ell+d} \cup
R^{small}_{\ell-d,\ell+d}  
\]
and put $\cYo=\cYo^5$.
\end{enumerate}
Having defined $\cYo$ we set 
\[
\cYjk = M\bigl( [k(\ell-d)+j ](\alpha-\beta),k(\ell+d+1)
\bigr) + \cYo\,,\;
(j,k)\in\cL \,.  
\]

\paragraph*{The percolation variables $U_{jk}$.}
Let $\cOjk$ be  the event that
for every parallelogram $\cP$ in $\cYjk$ there is an active path
in the graphical representation of 
the contact process which stays entirely in $\cP$ and
connects some point in the bottom edge of
$\cP$ to some point in the the top edge of $\cP$. 
Thus on $\cOjk$ there is some
point in the bottom edge of $\cYjk$ with the property that
there are active paths in $\cYjk$ connecting this point to
the top edge of every parallelogram in $\cYjk$, and in
particular to the top edges of the two top parallelograms 
$\cYjk$. This means that on $\cOjk$ there is a
point in the bottom edge of $\cYjk$ and active paths in
$\cYjk$ connecting this point to the bottom edges of
both $\cY_{j-1,k+1}$ and $\cY_{j+1,k+1}$. 

It is a consequence of Lemma~VI.3.17 in \cite{bible} that
$P(\O00)$ is 
close to 1 for large $M$. 
\begin{lemma}\label{lemma2}
For  $0 < \b < \a/3$ , $\lim_{M\to \infty}P(\cOo)=1$.
\end{lemma}

Proof: As in \cite{bible} let $\Ejk$ to be the event that
there is an active path in the graphical representation of the
contact process which goes from the bottom edge of $\Rjk$ to
the top edge, always staying entirely within $\Rjk$, and also
that there is an active path from the bottom edge of $\Ljk$ to
the top edge, always staying entirely within $\Ljk$.  It is
clear that the probability of connecting the bottom edge of a
small parallelogram to its top edge by an active path staying
in the parallelogram is bounded below by $P(\cE_{00})$. 
By Lemma 3.17 in \cite{bible}, for $0 < \b < \a/3$, 
$\lim_{M\to \infty}P(\cE_{00})=1$.  In the construction of
$\cYo$ there are most $h=2\ell+4d$ (if $d \geq1$) or $h=2\ell+1$ (if $d=0$)
parallelograms used. It 
follows from positive correlations that $P(\mathcal{O}_{00})
\geq 
P(\Ejk)^{h}$, and thus $\lim_{M\to
  \infty}P(\mathcal{O}_{00})=1$ \qed 

\bigskip
For $(j,k)\in \cL$ let $U_{jk}=1_{\cOjk}$. Then
$P(U_{jk}=1)=P(\mathcal{O}_{00})$ does not depend on $(j,k)$.
Furthermore, the $U_{jk}$ are 1-dependent, meaning that if 
$I\subset \cL$ is such that $\|(j,k)-(j',k')\|>1$ for all
$(j,k)\ne(j',k')\in I$, then the $U_{jk}, (j,k)\in I$ are
independent. This is because the corresponding space-time
regions $\cYjk,\cY_{j'k'}$ are disjoint. Using the $U_{jk}$ we
may construct a 1-dependent oriented percolation process in the
usual way. A {\em path} in $\cL$ is a sequence
$(j_1,k_1),...,(j_n,k_n)$ of points of $\cL$ which satisfies
$k_{i+1}=k_i+1$ and $j_{i+1}=j_i \pm 1$ for all $ 1 \leq i
\leq n-1$. The path is said to be {\em{open}} if
$U_{j_i,k_i}=1$ for each $ 1 \leq i \leq n-1$. It is clear
from the properties of the $\cOjk$ that if
$(j_1,k_1),...,(j_n,k_n)$ is an open path in $\cL$ then there
must an active path in the graphical representation from the
bottom edge of $\cY_{j_1,k_1}$ to the bottom edge of
$\cY_{j_n,k_n}$.

\begin{figure}[htp]
 \begin{center}
\includegraphics[totalheight=0.25\textheight]{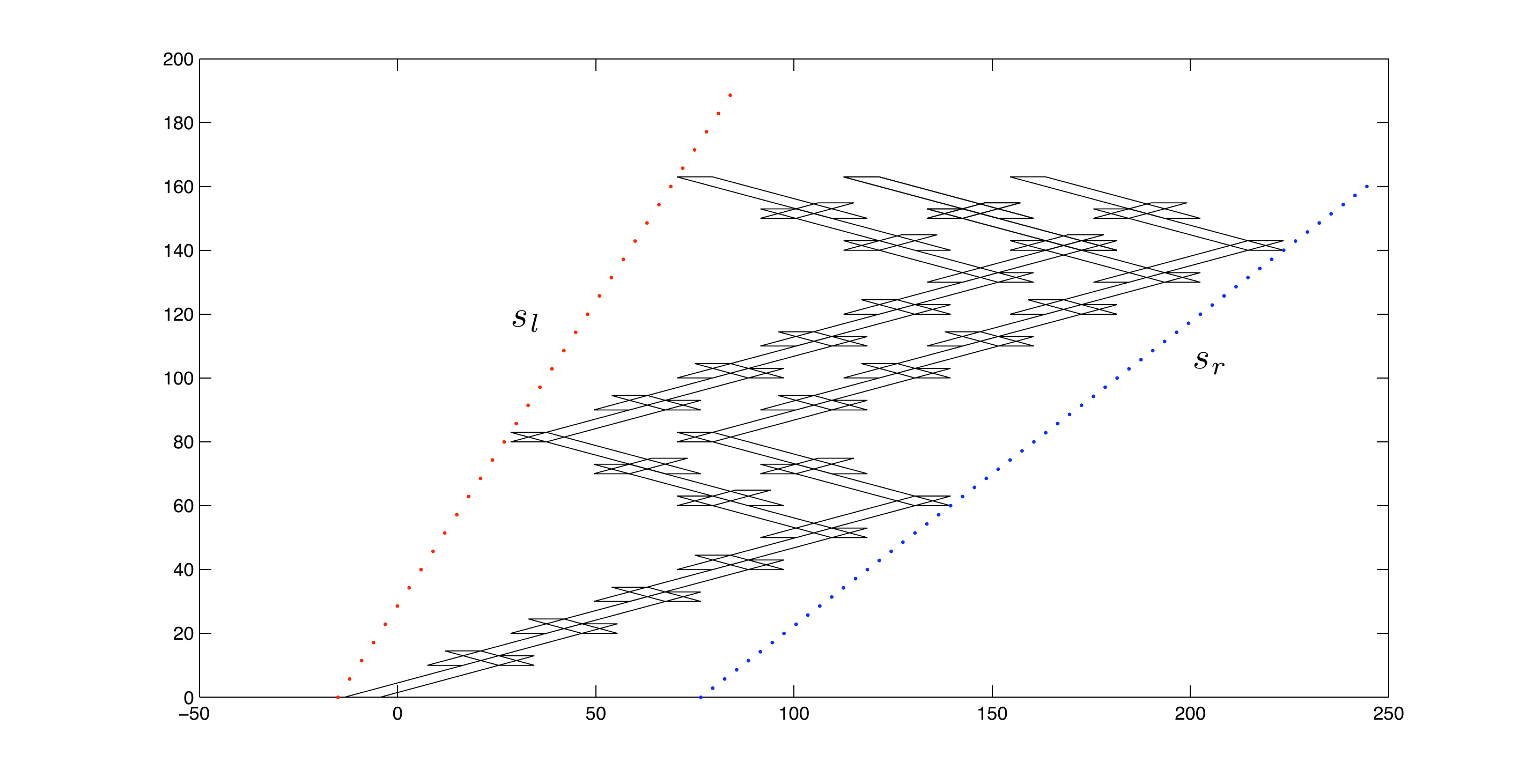}
 \caption{ $\cY_{00},\cY_{1,1},\cY_{-1,1}$}
 \end{center}
\end{figure}

If we let $\Omega_{\infty}$ be the event that there is an
infinite open path in $\cL$ starting at $(0,0)$, then by
Lemma~\ref{lemma2} above and Theorem~VI.3.19 of \cite{bible},
\begin{equation}\label{percolate}
\lim_{M\to\infty}P(\Omega_\infty)=1\;.
\end{equation}

\paragraph*{Survival of $\xi^\cW_t$.} Let 
$\cY =\cY(\ell,d,M) = \bigcup_{k=0}^\infty \bigcup_{j=-k}^k\cY_{jk}$.
On $\Omega_\infty$ there must be an infinite active
path in the graphical representation starting at some
$(x,0)$, $x\in[-3M\beta/2,-M\beta/2]$, which lies entirely
in $\cY$. Thus if $\cW$ is any space-time region such that
$\cY\subset\cW$, and $\xi^{\cW}_t$ is the $\cW$-restricted
contact process starting from $\{x:(x,0)\subset\cW\}$, then
$\xi^{\cW}_t\ne\emptyset \;\forall\;t\ge 0$ on 
$\Omega_\infty$. We will prove the following.

\medskip
\noindent{\bf Claim.} Assume \eqref{wedgespeeds} holds and
$\alpha=\alpha(\lambda)$. Then there exists $0<\beta<\alpha/3$
and integers $\ell',d'$ such that for all $M>0$, 
\begin{equation}\label{containment}
\cY(\ell',d',M/\alpha(\ell'+3))\subset
\cW(\alpha_l,\alpha_r,M) - 
(M/(\ell'+3),0) \;.
\end{equation}

\medskip\noindent
Given
\eqref{containment}, it follows from translation invariance
and \eqref{percolate} that 
\[
P(\xi^{\cW(\alpha_l,\alpha_r,M)}_t\ne\emptyset\;\forall\;t\ge 0) \ge
P(\Omega_\infty) \to 1 \text{ as }M\to\infty \,,
\]
proving \eqref{eqn:survival}.

To prove \eqref{containment} we first suppose that 
$\ell,d$,  are positive integers with $d<\ell$ and 
$M>0$. For
$(j,k) \in \cL$, the left upper corner of $\Ljk$ is
$(M(j(\a-\b)-\a -\b/2), M(k+1+\b/\a))$, and the right bottom
corner of $\Ljk$ is $(M(j(\a-\b)+3\b/2), Mk)$.  A little
thought shows that $\cY$ must be contained in the space-time
region bounded by the following two lines and the $x$-axis.
The first line connects the leftmost point of the top edge of $\cYo$ with
 the leftmost point of the top edge of $\cY_{-1,1}$, which are the
left upper corner of $L_{\ell-d,\ell+d}$ and the left upper
corner of $L_{2(\ell+d)-1,2(\ell+d)+1}$, namely, the points 
$(M((\ell-d)(\alpha-\beta)-\alpha-\beta/2,
M(\ell+d+1+\beta/\alpha)$ and $(M(2(\ell-d)(\a-\b)-2\a+\b/2),M(2(\ell+d+1)+\b/\a))$. 
The slope of this line is 
\begin{equation}\label{slope1}
s_l=\frac{\ell+d+1}{\ell-d-1}\frac{1}{\a-\b}
\end{equation} 
and it contains the point $(x_l,0)$ where $x_l= -M(3\beta /2 +
\beta/\alpha s_l)$.  The second line connects the
rightmost point of $\cYo$ with the rightmost point of
$\cY_{1,1}$, the bottom right corner of $L_{\ell+1,\ell+1}$
and the bottom right $L_{2(\ell+1)-d,2(\ell+1)+d}$, namely, the
points $(M((\ell+1)(\a-\b) +3 \b/2),M(\ell+1))$ and
$(M((2(\ell+1)-d)(\a-\b)+3\b/2),M(2(\ell+1)+d))$. The slope
of this line is
\begin{equation}\label{slope2}
s_r=\frac{\ell+d+1}{\ell-d+1}\frac{1}{\a-\b}
\end{equation} 
and it contains the point $(x_r,0)$ where $x_r=M((\ell+1)(\alpha-\beta-1/s_r)+3\beta/2)$. 

This analysis shows that $\cY(\ell,d,M)$ is contained in the
wedge $\cW(1/s_l,1/s_r, M')+(x_l,0)$, where $M' =
x_r-x_l$. A little algebra shows that 
$-M\alpha < x_l <x_r < M\alpha(\ell+2)$, 
and thus
\begin{equation}\label{containment2}
\cY(\ell,d,M)\subset
\cW(1/s_l,1/s_r,M\alpha(\ell+3))-(M\alpha ,0) \;.
\end{equation}
We now set $s_\ell=1/\alpha_\ell, s_r = 1/\alpha_r$ and 
solve \eqref{slope1} and \eqref{slope2} for $d$ and
$\ell$, obtaining 
\begin{equation}\label{dl}
\ell=\frac{s_r(s_l(\alpha-\beta)+1)}{s_l-s_r}\,,
\quad 
d=\frac{s_l(s_r(\alpha-\beta)-1)}{s_l-s_r}\;.
\end{equation}
Unfortunately, $\ell,d$ need not be integers. To deal with this problem
we first note that if $s_l \geq s_l' > s_r$ then for any
$M$, the wedge $\cW(\alpha_l,\alpha_r,M)$ contains the narrower
wedge $\cW(1/s'_{\ell},1/s_r,M)$. If we can find 
$s'_\ell$ and $0<\beta<\alpha/3$ such that 
\begin{equation}\label{intsols}
\ell'= \frac{s_r(s_l'(\alpha-\beta)+1)}{s_l' - s_r}
\text{ and }
d'= \frac{s_l'(s_r(\alpha-\beta)-1)}{s_l' - s_r}
\end{equation}
are both integers, then \eqref{containment} follows from 
\eqref{containment2}.

We can find $s'_\ell,\beta$ as follows. Let $m_0= 3/\alpha
s_r$ and take any integer $m>m_0$ such that $s_r
\frac{m}{m-1} < s_l$. Put $s_l' = s_r \frac{m}{m-1}$, so that 
$s_l>s'_l>s_r$.  Since
$m>3/\alpha s_r$, $ 1/3 \alpha m s_r >1$ and the interval
$( \frac{2}{3}\, \alpha\, m \, s_r, \alpha\, m\, s_r)$
must contain at least one integer. Since $\alpha s_r >1$, the
right endpoint of this interval is greater than $m$. Choose
any integer $c\ge m$ from the interval and 
put $\beta = \alpha - \frac{c}{m s_r}$.  Then 
$0< \beta < \alpha/3$ and $ s_r(\alpha -
\beta) = c/m$. A little algebra shows that $\ell',d'$ given in
\eqref{intsols} are the integers $\ell'=c+m-1, d'=c-m$, and we
are done.

\begin{figure}[htp]
\centering
\caption{ Wedge containing $\cY$ }
\includegraphics[totalheight=0.35\textheight]{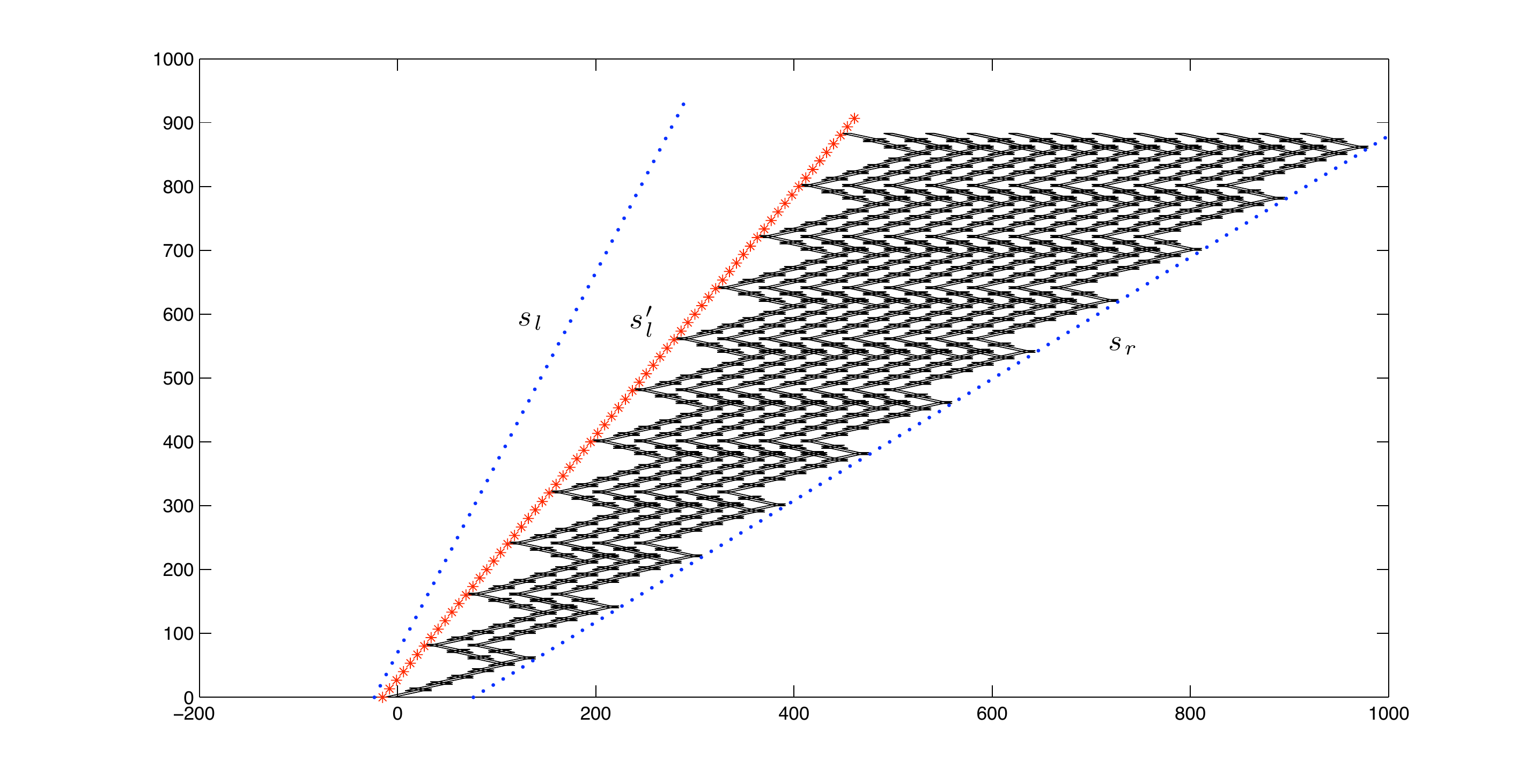}\label{wedge} 
\end{figure}

\section{Proof of Theorem~\ref{equilibrium}}
We begin by analyzing the
rightmost particle. Let $\cW(\alpha_r,M)=
\{(x,t) : t\ge 0, x\in (-\infty,M+\alpha_r t]\cap\Z\}$ 
and consider the restricted contact process
$\xi^{\cW(\alpha_r,M)}_t$ with initial state 
$\xi^{\cW(\alpha_r,M)}_0=(-\infty,M]\cap\Z$.
Let $\bar r_t$ be the right-edge process for $\xi^\cW_t$, 
$\bar r_t = \max\{x :
\xi^{\cW(\alpha_r,M)}_t(x) =1 \}$.
We claim that for every $M$, 
\begin{equation}\label{r-edge}
\lim_{t\to\infty} \dfrac{\bar r_t}t = 
\alpha_r ~ ~ a.s.
\end{equation}
By construction and \eqref{edgespeed},
$\limsup_{t\to\infty} \bar r_t/t\le 
\alpha_r$. For the lower bound,
fix $0<\vep<\alpha_r$ and define the 
region $\cW_\vep = \cW(\alpha_r-\vep,\alpha_r,M)$ and
restricted contact process $\xi^{\cW_\vep} _t$ with initial
state $\xi^{\cW_\vep} _0 = [0,M]\cap\Z$. Then
$\xi^{\cW_\vep}_t\subset \xi^{\cW(\alpha_r,M)}_t$, which
implies that on the event $\{\xi^{\cW_\vep}_t\ne \emptyset \;\forall\; t\ge 0\}$, 
$\liminf_{t\to\infty}\bar r_t/t \ge \alpha_r-\vep$. 
Theorem~\ref{survival} now implies we must have
$\liminf_{t\to\infty}\bar r_t/t \ge \alpha_r-\vep$ a.s.,
completing the proof of \eqref{r-edge}.  

It is a consequence of the nearest-neighbor interaction
mechanism that for any $\alpha_l<\alpha_r$ and $M$, with $\cW
= \cW(\alpha_l,\alpha_r,M)$,
\[
\xi^{\cW}_t(x) = \xi_t^{\cW(\alpha_r,M)}(x) \; \forall \; 
x\in[l^{\cW}_t,r^{\cW}_t] \text{ on }
\{\xi^{\cW}_t \ne\emptyset\}\,.
\]
This implies $r^{\cW}_t=\bar r_t$ on $\{\xi^{\cW}_t
\ne\emptyset\}$, and so by \eqref{r-edge}, 
$\lim_{t\to\infty}r^\cW_t/t =
\alpha_r$. We omit the similar
argument proving $\lim_{t\to\infty} l^\cW_t/t =
\alpha_l$.

For \eqref{nu-couple}, let $\xi^{\Z}_t$ denote the
unrestricted process with initial state $\xi^\Z_0=\Z$,
and let $\xi^\nu_t$ be the unrestricted process constructed as in
Section 2 with initial state
$\xi^\nu_0$ which has law $\nu$,  independent of the Poisson
processes. 
We observe again that the
nearest-neighbor interaction implies
\[
\xi^{\Z}_t (x) = 
\xi^{\cW}_t (x) \; \forall \; x\in [l^\cW_t,r^\cW_t]
\text{ on } \{\xi^\cW_t\ne\emptyset \; \forall \;t\ge 0\}\,.
\]
Standard exponential estimates for $P(\xi^\Z_t(x)\ne
\xi^\nu_t(x))=P(\xi^\Z_t(x)=1)-P(\xi^\nu_t(x)=1)$, a 
``filling in'' argument and Borel-Cantelli (see Theorem~I.2.30
of \cite{bible2}) 
imply that for any $A>0$,
\[
P( \xi^\Z_t = \xi^\nu_t \text{ on } [-At,At] \text{
  for all large } t )=1 
\]
Combining the above with \eqref{wedgegrowth} gives \eqref{nu-couple}. 

\section{Proof of Corollary~\ref{gbt}} 
We will make use of the graphical construction in
Section~\ref{harris} and define independent events
$\Omega_1,\Omega_2,\Omega_3$, 
each with positive probability, and such that
$\|\zeta_t\|_1\to\infty$ as $t\to\infty$ on their
intersection.  

First, since $\alpha(\lambda)$ is strictly increasing we may choose 
$\alpha(\lambda_2)<\alpha_l<\alpha_r<\alpha(\lambda_1)$. Fix
$M>2$ and write $\cW$ for $\cW(\alpha_l,\alpha_r,M)$. 
The first event is 
\[
\Omega_1=\{\text{there is no active 2-path from any $(x,0)$,
  $x<0$, to any point of $\cW(\alpha_l,\alpha_r,M)$}\} \;.
\]
Since the process of 2's is a contact process with
parameter $\lambda_2$,
and $\alpha(\lambda_2)<\alpha_l$, it
follows from \eqref{edgespeed} that
$\Omega_1$ has positive probability.

For the second event, choose $x_0\in\Z$ and $t_0>0$ such that 
$x_0=\alpha_lt_0$ and
$(x,t_0)\subset\cW$ for all $x\in
[x_0,x_0+M]\cap\Z$. Since $M>2$ the event,
\[
\Omega_2=\{\text{there is an active path in $\cW$
from $(0,0)$ to each of $(x,t_0), x\in[x_0,x_0+M]\cap\Z$}\}
\]
has positive probability. 

For the third event, define, for $t\ge t_0$,
\begin{multline*}
A_t=\{y: \text{there is an infinite active path in $\cW$
  from $(x,t_0)$ to $(y,t)$}\\
\text{for some $x\in[ x_0,x_0+M]\cap\Z$ }
\}
\end{multline*}
and put $\Omega_3=\{|A_t|\to\infty\text{ as }t\to\infty\}$. It
follows from Theorems~\ref{survival} and \ref{equilibrium}
that $\Omega_3$ has positive probability.

The events $\Omega_i$ are independent since they are defined in terms
of our Poisson processes over disjoint space-time regions.
Furthermore, it is easy to see from Remark~\ref{lingrowth}
that $\|\zeta_t\|_1\to\infty$ on their 
intersection, so we are done.

\vskip10mm
\noindent {\it
J. Theodore Cox \\
Department of Mathematics,
Syracuse University \hskip3mm
Syracuse, NY 13244\\ e-mail: {\tt jtcox@syr.edu}\\[3mm]
Nevena Mari\'{c} \\
Department of Mathematics and CS,
University of Missouri- St. Louis\hskip3mm St. Louis, MO 63121\\
e-mail: {\tt maricn@umsl.edu}\\[3mm]
Rinaldo Schinazi\\
Department of Mathematics, University of Colorado-Colorado Springs \hskip3mm
Colorado Springs, CO 80933\\
e-mail: {\tt rschinaz@uccs.edu}
}

\end{document}